\documentclass{article}
\usepackage{amsfonts}
\usepackage{amsmath}

\newtheorem{theorem}{Theorem}

\newtheorem{corollary}[theorem]{Corollary}

\newtheorem{lemma}{Lemma}

\begin{document}

\title{Critical branching processes evolving in an unfavorable random
environment\thanks{%
This work was supported by the Russian Science Foundation under grant
no.19-11-00111 https://rscf.ru/en/project/19-11-00111/ }}
\author{V.A.Vatutin\thanks{%
Steklov Mathematical Institute Gubkin street 8 119991 Moscow Russia Email:
vatutin@mi.ras.ru}, E.E.Dyakonova\thanks{%
Steklov Mathematical Institute Gubkin street 8 119991 Moscow Russia Email:
elena@mi-ras.ru}}
\maketitle

\begin{abstract}
Let $\left\{ Z_{n},n=0,1,2,...\right\} $ be a critical branching process in
random environment and let $\left\{ S_{n},n=0,1,2,...\right\} $ be its
associated random walk. It is known that if the increments of this random
walk belong (without centering) to the domain of attraction of a stable law,
then there exists a sequence $a_{1},a_{2},...,$ slowly varying at infinity
such that the conditional distributions
\begin{equation*}
\mathbf{P}\left( \frac{S_{n}}{a_{n}}\leq x\Big|Z_{n}>0\right) ,\quad x\in
(-\infty ,+\infty ),
\end{equation*}%
weakly converges, as $n\rightarrow \infty $ to the distribution of a
strictly positive and proper random variable. In this paper we supplement
this result with a description of the asymptotic behavior of the probability
\begin{equation*}
\mathbf{P}\left( S_{n}\leq \varphi (n);Z_{n}>0\right) ,
\end{equation*}%
if $\varphi (n)\rightarrow \infty $ \ as $n\rightarrow \infty $ in such a
way that $\varphi (n)=o(a_{n})$.

\end{abstract}

\section{Introduction and main result}

We consider critical branching processes evolving in an unfavorable random
environment. To formulate the problem under consideration and to give a
detailed description of our main result we start by recalling some basic
properties of such processes.

Let $\mathcal{M}$ be the space of all probability measures on $\mathbf{N}%
_{0}:=\{0,1,2,\ldots \}.$ Equipped with a metric, $\mathcal{M}$ becomes a
Polish space. Let $F$ be a random variable taking values in $\mathcal{M}$,
and let $F_{n},n\in \mathbf{N}:=\mathbf{N}_{0}\backslash \left\{ 0\right\} $
be a sequence of independent copies of $F$. The infinite sequence $\mathcal{E%
}=\left\{ F_{n},n\in \mathbf{N}\right\} $ is called a random environment.
Given the environment $\mathcal{E}$, we may construct the i.i.d. sequence of
generating functions
\begin{equation*}
F_{n}(s):=\sum_{j=0}^{\infty }F_{n}\left( \left\{ j\right\} \right)
s^{j},\quad s\in \lbrack 0,1].
\end{equation*}%
In the sequel we make no difference between an element $F_{n}\in \mathcal{M}$
and the respective generating function $F_{n}(s)$ and use for a random
variable $F,$ taking values in $\mathcal{M}$, the representation%
\begin{equation*}
F(s):=\sum_{j=0}^{\infty }F\left( \left\{ j\right\} \right) s^{j},\quad s\in
\lbrack 0,1].
\end{equation*}

A sequence of nonnegative random variables $\mathcal{Z}=\left\{ Z_{n},\ n\in
\mathbf{N}_{0}\right\} $ specified on\ a probability space $(\Omega ,%
\mathcal{F},\mathbf{P})$ is called a branching process in random environment
(BPRE), if $Z_{0}$ is independent of $\mathcal{E}$ and, given $\mathcal{E}$
the process $\mathbf{Z}$ is a Markov chain with
\begin{equation*}
\mathcal{L}\left( Z_{n}|Z_{n-1}=z_{n-1},\mathcal{E}=(f_{1},f_{2},...)\right)
=\mathcal{L}(\xi _{n1}+\ldots +\xi _{ny_{n-1}})
\end{equation*}%
for all $n\in \mathbf{N}$, $z_{n-1}\in \mathbf{N}_{0}$ and $%
f_{1},f_{2},...\in \mathcal{M}$, where $\xi _{n1},\xi _{n2},\ldots $ is a
sequence of i.i.d. random variables with distribution $f_{n}.$ Thus, $%
Z_{n-1} $ is the $(n-1)$th generation size of the population of the
branching process and $f_{n}$ is the offspring distribution of an individual
at generation $n-1$.

We denote $X_{i}=\log F_{i}^{\prime }(1),i=1,2,...$ and introduce the
sequence
\begin{equation*}
S_{0}=0,\quad S_{n}=X_{1}+...+X_{n},\ n\geq 1,
\end{equation*}%
is called the associated random walk for the process $\mathcal{Z}$.

We need the subset
\begin{equation*}
\mathcal{A}=\{0<\alpha <1;\,|\beta |<1\}\cup \{1<\alpha <2;|\beta |\leq
1\}\cup \{\alpha =1,\beta =0\}\cup \{\alpha =2,\beta =0\}
\end{equation*}%
of the set $\mathbb{R}^{2}.$ For a pair $(\alpha ,\beta )\in \mathcal{A}$
and a random variable $X$ we write $X\in \mathcal{D}\left( \alpha ,\beta
\right) $ if the distribution of $X$ belongs (without centering) to the
domain of attraction of a stable law with density $g_{\alpha .\beta
}(x),x\in (-\infty ,+\infty )$ and the characteristic function%
\begin{equation*}
G_{\alpha ,\beta }\mathbb{(}w\mathbb{)}=\int_{-\infty }^{+\infty
}e^{iwx}g_{\alpha .\beta }(x)\,dx=\exp \left\{ -c|w|^{\,\alpha }\left(
1-i\beta \frac{w}{|w|}\tan \frac{\pi \alpha }{2}\right) \right\} ,\ c>0,
\end{equation*}%
This implies, in particular, that there is an increasing sequence of
positive numbers
\begin{equation}
a_{n}\ =\ n^{1/\alpha }\ell (n),  \label{defA}
\end{equation}%
where $\ell (1),\ell (2),\ldots $ is a slowly varying sequence, such that,
as $n\rightarrow \infty $
\begin{equation*}
\mathcal{L}\left\{ \frac{S_{nt}}{a_{n}},t\geq 0\right\} \overset{D}{%
\rightarrow }\mathcal{L}\left\{ Y_{t},t\geq 0\right\} ,
\end{equation*}%
where the symbol \ \ $\overset{D}{\rightarrow }$ stands for the convergence
in distribution in the space $D[0,+\infty )$ with Skorokhod topology and the
process $\mathcal{Y}=\left\{ Y_{t},t\geq 0\right\} $ is strictly stable and
has marginal distributions specified by the characteristic functions \ \ \ \
\ \ \ \ \ \ \
\begin{equation*}
\mathbf{E}e^{iwY_{t}}=G_{\alpha ,\beta }\mathbb{(}wt^{1/\alpha }\mathbb{)}%
,t\geq 0.
\end{equation*}%
Observe that if $X_{n}\overset{d}{=}X\in \mathcal{D}\left( \alpha ,\beta
\right) $ for all $n\in \mathbf{N}$ then (see,$\ $\cite{Zol57} or \cite[p.
380]{BGT87}) the limit%
\begin{equation*}
\lim_{n\rightarrow \infty }\mathbf{P}\left( S_{n}>0\right) =\rho =\mathbf{P}%
\left( Y_{1}>0\right)
\end{equation*}%
exists, where
\begin{equation*}
\displaystyle\rho =\frac{1}{2}+\frac{1}{\pi \alpha }\arctan \left( \beta
\tan \frac{\pi \alpha }{2}\right) .
\end{equation*}

We now formulate our first restriction on the properties of the BPRE.

\paragraph{Condition B1.}

\emph{The random variables }$X_{n}=\log F_{n}^{\prime }(1),n\in \mathbf{N}$%
\emph{\ are independent copies of a random variable }$X\in \mathcal{D}\left(
\alpha ,\beta \right) $\emph{,} $\left\vert \beta \right\vert <1,$ \emph{%
whose distribution is absolutely continuous. }

Our second assumption on the environment concerns the standardized truncated
second moment of the generating function $F$:
\begin{equation*}
\zeta (b)=\frac{\sum_{k=b}^{\infty }k^{2}F\left( \left\{ k\right\} \right) }{%
\left( \sum_{i=b}^{\infty }iF\left( \left\{ i\right\} \right) \right) ^{2}}.
\end{equation*}

\paragraph{Condition B2.}

\emph{There exist $\varepsilon >0$ and $b\in $}$\mathbf{N}$ \emph{such that}
\begin{equation*}
\mathbf{E}[(\log ^{+}\zeta (b))^{\alpha +\varepsilon }]\ <\ \infty \ ,
\end{equation*}%
\emph{where }$\log ^{+}x=\log (\max (x,1))$\emph{.}

It is known (see \cite[Theorem 1.1 and Corollary 1.2]{agkv}) that if
Conditions B1, B2 are valid then there exist a number $\theta \in (0,\infty
) $ and a sequence $l(1),l(2)...$ slowly varying at infinity such that, as $%
n\rightarrow \infty $%
\begin{equation*}
\mathbf{P}\left( Z_{n}>0\right) \sim \theta \mathbf{P}\left( \min \left(
S_{1},...,S_{n}\right) \geq 0\right) \sim \theta n^{-(1-\rho )}l(n).
\end{equation*}%
Besides (see \cite[Theorem 1.5]{agkv}), for any $t\in \left[ 0,1\right] $
and any $x\geq 0$
\begin{equation*}
\lim_{n\rightarrow \infty }\mathbf{P}\left( \frac{S_{nt}}{a_{n}}\leq
x|Z_{n}>0\right) =\mathbf{P}\left( Y^{+}(t)\leq x\right) ,
\end{equation*}%
where $\mathcal{Y}^{+}=\left\{ Y^{+}(t),0\leq t\leq 1\right\} $ denotes the
meander of the strictly stable process $\mathcal{Y}$.

Thus, the associated random walk that provides survival of the population to
a distant moment $n$ \ grows like $a_{n}$ times a random positive multiplier.

Our aim is to investigate the behavior of the probability of the event $%
\left\{ Z_{n}>0\right\} $ in an unfavorable environment, namely, when $0\leq
S_{n}=o(a_{n})$ as $n\rightarrow \infty $.

To formulate the desired statement we denote%
\begin{equation*}
M_{n}=\max \left( S_{1},...,S_{n}\right) ,\quad L_{n}=\min \left(
S_{1},...,S_{n}\right) ,
\end{equation*}%
and introduce right-continuous renewal functions

\begin{align*}
U(x)& =I\left\{ x\geq 0\right\} +\sum_{n=1}^{\infty }\mathbf{P}\left(
S_{n}\geq -x,M_{n}<0\right) ,\ x\in \mathbb{R}, \\
V(x)\ & =\ I\left\{ x<0\right\} +\sum_{k=1}^{\infty }\mathbf{P}\left(
S_{n}<-x,L_{n}\geq 0\right) ,\ x\in \mathbb{R}.
\end{align*}%
It is well-known that $U(x)=O(x)$ and $V(-x)=O(x)$ as $x\rightarrow \infty $.

\textbf{Remark 1.} Observe that if Condition B1 is valid then
\begin{equation*}
V(x)\ =\tilde{V}(x)\ =\ I\left\{ x<0\right\} +\sum_{k=1}^{\infty }\mathbf{P}%
\left( S_{n}<-x,L_{n}>0\right) ,\ x\in \mathbb{R}.
\end{equation*}

We use this fact many times below referring to the results of \cite{VW09}.

We now formulate the main result of the note.

\begin{theorem}
\label{T_smallDevi}Let Conditions B1, B2 be valid. If a function $\varphi
(n)\rightarrow \infty $ as $n\rightarrow \infty $ and $\varphi (n)=o(a_{n})$
then%
\begin{equation*}
\mathbf{P}\left( Z_{n}>0;S_{n}\leq \varphi (n)\right) \sim \Theta \frac{%
g_{\alpha ,\beta }(0)}{na_{n}}\int_{0}^{\varphi (n)}V(-z)dz
\end{equation*}%
as $n\rightarrow \infty $, where $\Theta $ is a positive constant specified
in formula (\ref{Def_Theta}) below.
\end{theorem}

Theorem \ref{T_smallDevi} compliments Theorem 1.1 in \cite{agkv} where the
asymptotic behavior of the survival probability $\mathbf{P}\left(
Z_{n}>0\right) $ was investigated as $n\rightarrow \infty $.

In the sequel we denote by $C,C_{1},C_{2},...,$ some positive constants that
do not necessarily coincide in different formulas.

\section{Auxiliary results}

Proving Theorem \ref{T_smallDevi} we will use random walks that start from
any point $x\in \mathbb{R}$. In such cases we write the respective
probabilities as $\mathbf{P}_{x}\left( \cdot \right) $. We also write $%
\mathbf{P}$ instead of $\mathbf{P}_{0}$.

We now formulate a number of statements that show importance of the
functions $U$ and $V$.

\begin{lemma}
\label{L_estimV} If Condition B1 is valid than for any $\lambda \in
(0,\infty )$ there exists a constant $C(\lambda )$ such that
\begin{equation*}
\int_{0}^{\lambda \varphi (n)}V(-z)dz\leq C(\lambda )\int_{0}^{\varphi
(n)}V(-z)dz
\end{equation*}%
for all $n\geq 1$.
\end{lemma}

\textbf{Proof}. If $X\in \mathcal{D}(\alpha ,\beta ),$ then (compare with
Lemma 13 in \cite{VW09})
\begin{equation*}
V(-x)=x^{\alpha \rho }l_{2}(x)
\end{equation*}%
as $x\rightarrow \infty $, where $l_{2}(x)$ is a slowly varying function.
Therefore, the function
\begin{equation*}
\int_{0}^{x}V(-z)dz
\end{equation*}%
is regularly varying of the index $\alpha \rho +1$ as $x\rightarrow \infty $
(see \cite[Ch. VIII, Sec. 9, Theorem 1]{Fel2008}). Now the statement of the
lemma follows from the properties of regularly varying \ functions with a
positive index.

The lemma is proved.

Let
\begin{equation}
b_{n}=\frac{1}{a_{n}n}=\frac{1}{\ n^{1/\alpha +1}\ell (n)}.  \label{Def_b}
\end{equation}

\begin{lemma}
\label{L_double} If the distribution of the random \ variable $X$ satisfies
Condition B1 then there is a constant $C>0$ such that for all $n\geq 1$ and
for all $x,y\geq 0$
\begin{equation}
\mathbf{P}_{x}\left( 0\leq S_{n}<y,L_{n}\geq 0\right) \ \leq \
C\,b_{n}\,U(x)\int_{0}^{y}V(-z)dz\ ,  \label{Rough1}
\end{equation}%
and for all $x,y\leq 0$
\begin{equation}
\mathbf{P}_{x}\left( y\leq S_{n}<0,M_{n}<0\right) \ \leq \
C\,b_{n}\,V(x)\int_{y}^{0}U(-z)\ dz.  \label{Rough2}
\end{equation}
\end{lemma}

\textbf{Proof}. According to Proposition 2.3 in \cite{ABGV2011} there is a
constant $C>0$ such that for all $n\geq 1$ and for all $x,z\geq 0$
\begin{equation*}
\mathbf{P}_{x}\left( z-1\leq S_{n}<z,L_{n}\geq 0\right) \ \leq \
C\,b_{n}\,U(x)V(-z)\ ,
\end{equation*}%
and for all $x,z\leq 0$
\begin{equation*}
\mathbf{P}_{x}\left( z\leq S_{n}<z+1,M_{n}<0\right) \ \leq \
C\,b_{n}\,V(x)U(-z)\ .
\end{equation*}

Integration with respect to $z$ of the first inequality over the interval $%
[0,y),y\geq 0,$ and the second inequality over the interval $(y.0],y<0$
gives the desired statement.

The lemma is proved.

The next theorem is a restatement of Theorem 4 in \cite{VW09} in our
notation and refines (\ref{Rough1}) for $x=0$.

\begin{theorem}
\label{T_VatWach} If the distribution of the random \ variable $X$ satisfies
Condition B1. Then for any $\Delta >0$%
\begin{equation*}
\mathbf{P}\left( S_{n}\in \lbrack y,y+\Delta \right) ,L_{n}\geq 0)\sim
g_{\alpha ,\beta }(0)b_{n}\int_{y}^{y+\Delta }V(-w)dw
\end{equation*}%
uniformly in \thinspace $y\in (0,\delta _{n}a_{n}]$ where $\delta
_{n}\rightarrow 0$ as $n\rightarrow \infty $.
\end{theorem}

Integration over $y\in (0,x]$ the relation shown in Theorem \ref{T_VatWach}
leads to the following important conclusion.

\begin{corollary}
\label{C_IntegVW} Under the conditions of Theorem \ref{T_VatWach}%
\begin{equation*}
\mathbf{P(}S_{n}\leq x,L_{n}\geq 0)\sim g_{\alpha ,\beta
}(0)b_{n}\int_{0}^{x}V(-w)dw
\end{equation*}%
uniformly in \thinspace $x\in (0,\delta _{n}a_{n}],$ where $\delta
_{n}\rightarrow 0$ as $n\rightarrow \infty $.
\end{corollary}

We need some refinements of Theorem \ref{T_VatWach} established by Doney.
For $x\geq 0$ and $y\geq 0$ we write $\ $\ for brevity $x_{n}$ for $x/a_{n}$
and $y_{n}$ for $y/a_{n}$.

\begin{lemma}
\label{DonLocal} (\cite[Proposition 18]{Don12}) If the distribution of the
random \ variable $X$ satisfies Condition B1, then for any $\Delta >0$
\begin{equation*}
\mathbf{P}_{x}\left( S_{n}\in \lbrack y,y+\Delta ),L_{n}\geq 0\right) \sim
g_{\alpha ,\beta }(0)b_{n}U(x)\int_{y}^{y+\Delta }V(-w)dw
\end{equation*}%
uniformly with respect to $x,y\geq 0$ such that $\max
(x_{n},y_{n})\rightarrow 0$ as $n\rightarrow \infty .$
\end{lemma}

In view of the identities
\begin{equation*}
\begin{array}{rl}
\mathbf{E}[U(x+X);X+x\geq 0]\ =\ U(x)\ , & x\geq 0\ , \\
\mathbf{E}[V(x+X);X+x<0]\ =\ V(x)\ , & x\leq 0\ ,%
\end{array}%
\end{equation*}%
which hold for any oscillating random walk (see \cite[Ch. 4.4.3]{KV2017}) $U$
and $V$ give rise to two new probability measures $\mathbf{P}^{+}$ and $%
\mathbf{P}^{-}$. The construction procedure is standard and explained for $%
\mathbf{P}^{+}$ and $\mathbf{P}^{-}$ in detail in \cite{agkv} and \cite%
{ABGV2011}, respectively (see also \cite[Ch. 5.2]{KV2017}). We recall here
only some basic definitions related to the construction.

We assume that the random walk $\mathcal{S}=\left\{ S_{n},n\geq 0\right\} $
is adapted to some filtration $\mathcal{F}=(\mathcal{F}_{n},n\geq 0)$ and
construct probability measures $\mathbf{P}_{x}^{+}$, $x\geq 0,$ satisfying
for any bounded and measurable function $g:\mathcal{S}^{n+1}\rightarrow
\mathbb{R}$ the equality%
\begin{equation*}
\mathbf{E}_{x}^{+}[g(R_{0},\ldots ,R_{n})]\ =\ \frac{1}{U(x)}\mathbf{E}%
_{x}[g(R_{0},\ldots ,R_{n})U(S_{n});L_{n}\geq 0]\ ,\ n\in \mathbf{N}_{0},
\end{equation*}%
where $R_{0},R_{1},\ldots $ is a sequence of $\mathcal{S}$-valued random
variables, adapted to the filtration $\mathcal{F}$.

Similarly $V$ gives rise to probability measures $\mathbf{P}_{x}^{-}$, $%
x\leq 0$, characterized for each $n\in \mathbf{N}_{0}$ by the equation
\begin{equation*}
\mathbf{E}_{x}^{-}[g(R_{0},\ldots ,R_{n})]\ :=\ \frac{1}{V(x)}\mathbf{E}%
_{x}[g(R_{0},\ldots ,R_{n})V(S_{n});M_{n}<0]\ .
\end{equation*}

Now we prove the following statement which generalizes Lemma 2.5 in \cite%
{agkv}.

\begin{lemma}
\label{L_cond} Assume Condition $B1$. Let $H_{1},H_{2},...,$ be a uniformly
bounded sequence of random variables adapted to the filtration $\mathcal{%
\tilde{F}=}\left\{ \mathcal{\tilde{F}}_{k},k\in \mathbf{N}\right\} $, which
converges $\mathbf{P}^{+}$-a.s. to a random variable $H_{\infty }$ as $%
n\rightarrow \infty $. Suppose that $\varphi (n),$ $n\in \mathbf{N}$ is a
function such that $\inf_{n\in \mathbf{N}}\varphi (n)\geq C>0$ and $\varphi
(n)=o(a_{n})$ as $n\rightarrow \infty $. Then%
\begin{equation*}
\lim_{n\rightarrow \infty }\frac{\mathbf{E}\left[ H_{n};S_{n}\leq \varphi
(n),L_{n}\geq 0\right] }{\mathbf{P}\left( S_{n}\leq \varphi (n),L_{n}\geq
0\right) }=\mathbf{E}^{+}\left[ H_{\infty }\right] .
\end{equation*}
\end{lemma}

\textbf{Proof}. For a fixed $k<n$ we have
\begin{equation*}
\frac{\mathbf{E}\left[ H_{k};S_{n}\leq \varphi (n),L_{n}\geq 0\right] }{%
\mathbf{P}\left( S_{n}\leq \varphi (n),L_{n}\geq 0\right) }=\mathbf{E}\left[
H_{k}\frac{\mathbf{P}_{S_{k}}(S_{n-k}^{\prime }\leq \varphi
(n),L_{n-k}^{\prime }\geq 0)}{\mathbf{P}\left( S_{n}\leq \varphi
(n),L_{n}\geq 0\right) };L_{k}\geq 0\right] ,
\end{equation*}%
where $\mathcal{S}^{\prime }=\left\{ S_{n}^{\prime },n=0,1,2,...\right\} $
is a random walk being a probabilistic copy of the random walk $\mathcal{S}$
and is independent of the set$\mathcal{\ }\left\{ S_{j},j=0,1,...,k\right\} $%
. We know by Lemma \ref{L_double}, Corollary \ref{C_IntegVW} and the
condition $\inf_{n\in \mathbf{N}}\varphi (n)\geq C>0$ that, there exist
constants $C,C_{1}$ and $C_{2}$ such that for any fixed $k$ and all $n\geq k$
\begin{equation*}
\frac{\mathbf{P}_{S_{k}}(S_{n-k}^{\prime }\leq \varphi (n),L_{n-k}^{\prime
}\geq 0)}{\mathbf{P}\left( S_{n}\leq \varphi (n),L_{n}\geq 0\right) }\leq
\frac{C\,b_{n-k}\,U(S_{k})\int_{0}^{\varphi (n)}V(-z)dz}{C_{1}\,b_{n}\,%
\int_{0}^{\varphi (n)}V(-z)dz}\leq C_{2}U(S_{k}).
\end{equation*}%
Further, recalling Corollary\textbf{\ \ref{C_IntegVW} }and the definition (%
\ref{Def_b}) we see that, for each fixed $x\geq 0$ and $k\in \mathbf{N}$
\begin{equation*}
\lim_{n\rightarrow \infty }\frac{\mathbf{P}_{x}(S_{n-k}^{\prime }\leq
\varphi (n),L_{n-k}^{\prime }\geq 0)}{\mathbf{P}\left( S_{n}\leq \varphi
(n),L_{n}\geq 0\right) }=\lim_{n\rightarrow \infty }\frac{%
\,b_{n-k}\,U(x)\int_{0}^{\varphi (n)}V(-z)dz}{b_{n}\,\int_{0}^{\varphi
(n)}V(-z)dz}=U(x).\text{ }
\end{equation*}%
Since
\begin{equation*}
\mathbf{E}\left[ H_{k}U(S_{k});L_{k}\geq 0\right] =\mathbf{E}^{+}\left[ H_{k}%
\right] <\infty ,
\end{equation*}%
it follows by the dominated convergence theorem that, for each fixed $k$%
\begin{eqnarray*}
&&\lim_{n\rightarrow \infty }\frac{\mathbf{E}\left[ H_{k};S_{n}\leq \varphi
(n),L_{n}\geq 0\right] }{\mathbf{P}\left( S_{n}\leq \varphi (n),L_{n}\geq
0\right) } \\
&=&\mathbf{E}\left[ H_{k}\times \lim_{n\rightarrow \infty }\frac{\mathbf{P}%
_{S_{k}}(S_{n-k}^{\prime }\leq \varphi (n),L_{n-k}^{\prime }\geq 0)}{\mathbf{%
P}\left( S_{n}\leq \varphi (n),L_{n}\geq 0\right) };L_{k}\geq 0\right] \\
&=&\mathbf{E}\left[ H_{k}U(S_{k});L_{k}\geq 0\right] =\mathbf{E}^{+}\left[
H_{k}\right] .
\end{eqnarray*}%
Further, in view of the estimate (\ref{Rough1}) for each $\lambda >1$ we
have
\begin{eqnarray*}
&&\left\vert \mathbf{E}\left[ \left( H_{n}-H_{k}\right) ;S_{\lambda n}\leq
\varphi (n),L_{\lambda n}\geq 0\right] \right\vert \\
&\leq &\mathbf{E}\left[ \left\vert H_{n}-H_{k}\right\vert \mathbf{P}%
_{S_{n}}(S_{n(\lambda -1)}^{\prime }\leq \varphi (n),L_{n(\lambda
-1)}^{\prime }\geq 0);L_{n}\geq 0\right] \\
&\leq &Cb_{n(\lambda -1)}\int_{0}^{\varphi (n)}V(-z)dz\times \mathbf{E}\left[
\left\vert H_{n}-H_{k}\right\vert \,\,U(S_{n}),L_{n}\geq 0\right] \\
&=&Cb_{n(\lambda -1)}\int_{0}^{\varphi (n)}V(-z)dz\times \mathbf{E}^{+}\left[
\left\vert H_{n}-H_{k}\right\vert \,\right] .
\end{eqnarray*}%
Hence, using (\ref{Def_b}) and Corollary \ref{C_IntegVW} we conclude that%
\begin{eqnarray*}
\left\vert \mathbf{E}\left[ (H_{n}-H_{k})|S_{\lambda n}\leq \varphi
(n),L_{\lambda n}\geq 0\right] \right\vert &\leq &C\mathbf{E}^{+}\left[
\left\vert H_{n}-H_{k}\right\vert \right] \frac{b_{n(\lambda
-1)}\int_{0}^{\varphi (n)}V(-z)dz}{C_{1}\,b_{n\lambda }\,\int_{0}^{\varphi
(n)}V(-z)dz} \\
&\leq &C_{2}\left( \frac{\lambda }{\lambda -1}\right) ^{1+1/\alpha }\mathbf{E%
}^{+}\left[ \left\vert H_{n}-H_{k}\right\vert \right] .
\end{eqnarray*}%
Letting first $n$ and then $k$ to infinity we see that for each $\lambda >1$
the right-hand side of the previous relation vanishes by the dominated
convergence theorem.

Using this result we see that%
\begin{eqnarray*}
\lim_{n\rightarrow \infty }\mathbf{E}\left[ H_{n}|S_{\lambda n}\leq \varphi
(n),L_{\lambda n}\geq 0\right] &=&\lim_{k\rightarrow \infty
}\lim_{n\rightarrow \infty }\frac{\mathbf{E}\left[ \left( H_{n}-H_{k}\right)
;S_{\lambda n}\leq \varphi (n),L_{\lambda n}\geq 0\right] }{\mathbf{P}\left(
S_{\lambda n}\leq \varphi (n),L_{\lambda n}\geq 0\right) } \\
&&+\lim_{k\rightarrow \infty }\lim_{n\rightarrow \infty }\frac{\mathbf{E}%
\left[ H_{k};S_{\lambda n}\leq \varphi (n),L_{\lambda n}\geq 0\right] }{%
\mathbf{P}\left( S_{\lambda n}\leq \varphi (n),L_{\lambda n}\geq 0\right) }
\\
&=&\lim_{k\rightarrow \infty }\mathbf{E}^{+}\left[ H_{k}\right] =\mathbf{E}%
^{+}\left[ H_{\infty }\right] ,
\end{eqnarray*}%
which we rewrite as
\begin{equation*}
\mathbf{E}\left[ H_{n};S_{\lambda n}\leq \varphi (n),L_{\lambda n}\geq 0%
\right] =\left( \mathbf{E}^{+}\left[ H_{\infty }\right] +o(1)\right) \mathbf{%
P}\left( S_{\lambda n}\leq \varphi (n),L_{\lambda n}\geq 0\right) .
\end{equation*}%
Asuming without loss of generality that \ $\mathbf{E}^{+}\left[ H_{\infty }%
\right] \leq 1$ we conclude that
\begin{eqnarray*}
&&\left\vert \mathbf{E}\left[ H_{n};S_{n}\leq \varphi (n),L_{n}\geq 0\right]
-\mathbf{E}^{+}\left[ H_{\infty }\right] \mathbf{P}\left( S_{\lambda n}\leq
\varphi (n),L_{\lambda n}\geq 0\right) \right\vert \\
&\leq &\left\vert \mathbf{E}\left[ H_{n};S_{\lambda n}\leq \varphi
(n),L_{\lambda n}\geq 0\right] -\mathbf{E}^{+}\left[ H_{\infty }\right]
\mathbf{P}\left( S_{\lambda n}\leq \varphi (n),L_{\lambda n}\geq 0\right)
\right\vert \\
&&+\left\vert \mathbf{P}\left( S_{\lambda n}\leq \varphi (n),L_{\lambda
n}\geq 0\right) -\mathbf{P}\left( S_{n}\leq \varphi (n),L_{n}\geq 0\right)
\right\vert .
\end{eqnarray*}%
We have proved that the first summand at the right-hand side of the
inequality is of order $o\left( \mathbf{P}\left( S_{\lambda n}\leq \varphi
(n),L_{\lambda n}\geq 0\right) \right) $ as $n\rightarrow \infty $, and,
therefore, the order $o\left( \mathbf{P}\left( S_{n}\leq \varphi
(n),L_{n}\geq 0\right) \right) ,$ since
\begin{equation*}
\lim_{n\rightarrow \infty }\frac{\mathbf{P}\left( S_{\lambda n}\leq \varphi
(n),L_{\lambda n}\geq 0\right) }{\mathbf{P}\left( S_{n}\leq \varphi
(n),L_{n}\geq 0\right) }=C(\lambda )<\infty
\end{equation*}%
by Corollary \ref{C_IntegVW} and Lemma \ref{L_estimV}.

Further, again by Corollary \ref{C_IntegVW} and the definition (\ref{Def_b})
we have
\begin{eqnarray*}
&&\left\vert \mathbf{P}\left( S_{\lambda n}\leq \varphi (n),L_{\lambda
n}\geq 0\right) -\mathbf{P}\left( S_{n}\leq \varphi (n),L_{n}\geq 0\right)
\right\vert \\
&\leq &\left\vert \mathbf{P}\left( S_{\lambda n}\leq \varphi (n),L_{\lambda
n}\geq 0\right) -g_{\alpha ,\beta }(0)b_{n\lambda }\int_{0}^{\varphi
(n)}V(-z)dz\right\vert \\
&&+\left\vert \mathbf{P}\left( S_{n}\leq \varphi (n),L_{n}\geq 0\right)
-g_{\alpha ,\beta }(0)b_{n}\int_{0}^{\varphi (n)}V(-z)dz\right\vert \\
&&+g_{\alpha ,\beta }(0)\left\vert b_{n\lambda }-b_{n}\right\vert
\int_{0}^{\varphi (n)}V(-z)dz \\
&=&o\left( b_{n}\int_{0}^{\varphi (n)}V(-z)dz\right) +g_{\alpha ,\beta
}(0)b_{n}\left\vert \frac{b_{n\lambda }}{b_{n}}-1\right\vert
\int_{0}^{\varphi (n)}V(-z)dz.
\end{eqnarray*}%
Hence, letting $\lambda \downarrow 1$ we see that%
\begin{equation*}
\lim_{\lambda \downarrow 1}\lim_{n\rightarrow \infty }\frac{\left\vert
\mathbf{P}\left( S_{\lambda n}\leq \varphi (n),L_{\lambda n}\geq 0\right) -%
\mathbf{P}\left( S_{n}\leq \varphi (n),L_{n}\geq 0\right) \right\vert }{%
b_{n}\int_{0}^{\varphi (n)}V(-z)dz}=0.
\end{equation*}%
Combining the obtained estimates we get the statement of the lemma.

\section{Proof of Theorem \protect\ref{T_smallDevi}}

Introduce iterations of probability generating functions $\
F_{1}(.),F_{2}(.),...,$ by setting
\begin{equation*}
F_{k,n}(s)=F_{k+1}(F_{k+2}(\ldots (F_{n}(s))\ldots ))
\end{equation*}%
for $0\leq k\leq n-1$, $0\leq s\leq 1$, and letting $F_{n,n}(s)=s.$ Using
this notation we write%
\begin{equation*}
\mathbf{P}\left( Z_{n}>0|\ F_{k+1},\ldots ,F_{n};Z_{k}=1\right)
=1-F_{k,n}(0).
\end{equation*}%
In particular, for any $j\leq n$%
\begin{eqnarray*}
1-F_{0,n}(0) &=&\mathbf{P}\left( Z_{n}>0|\ F_{1},\ldots ,F_{n};Z_{0}=1\right)
\\
&\leq &\mathbf{P}\left( Z_{j}>0|\ F_{1},\ldots ,F_{j};Z_{0}=1\right) \\
&=&1-F_{0,j}(0)\leq e^{S_{j}}.
\end{eqnarray*}%
Therefore, $\lim_{n\rightarrow \infty }F_{0,n}(0)=F_{0,\infty }(0)$ exists
a.s.

We set%
\begin{equation*}
\tau _{n}=\min \left\{ 0\leq k\leq n:S_{k}=\min (0,L_{n})\right\}
\end{equation*}%
and write for $1<J<n$ the representation%
\begin{equation*}
\mathbf{P}\left( Z_{n}>0;S_{n}\leq \varphi (n)\right) =\sum_{j=0}^{J}\mathbf{%
P}\left( Z_{n}>0;S_{n}\leq \varphi (n),\tau _{n}=j\right) +R(J,n),
\end{equation*}%
where%
\begin{equation*}
R(J,n)=\sum_{j=J+1}^{n}\mathbf{P}\left( Z_{n}>0;S_{n}\leq \varphi (n),\tau
_{n}=j\right) .
\end{equation*}

\begin{lemma}
\label{L_remainder} If $\varphi (n)\rightarrow \infty $ then under the
condition of Theorem \ref{T_smallDevi}
\begin{equation*}
\lim_{J\rightarrow \infty }\lim_{n\rightarrow \infty }\frac{R(J,n)}{\mathbf{P%
}\left( S_{n}\leq \varphi (n),L_{n}\geq 0\right) }=0.
\end{equation*}
\end{lemma}

\textbf{Proof}. For each $j\in \mathbf{N}_{0}$ we have%
\begin{eqnarray*}
&&\mathbf{P}\left( Z_{n}>0;S_{n}\leq \varphi (n),\tau _{n}=j\right) =\mathbf{%
E}\left[ 1-F_{0,n}(0);S_{n}\leq \varphi (n),\tau _{n}=j\right] \\
&\leq &\mathbf{E}\left[ e^{S_{j}};S_{n}\leq \varphi (n);\tau _{n}=j\right] .
\end{eqnarray*}%
By the duality property of random walks and relations (\ref{Rough1}) and (%
\ref{Rough2})\ we get for $k>0$ and $y>0$ the inequalities%
\begin{equation}
\mathbf{P}\left( S_{n}\in \lbrack -k,-k+1),\tau _{n}=n\right) =\mathbf{P}%
\left( S_{n}\in \lbrack -k,-k+1),M_{n}<0\right) \leq C\,b_{n}\,U(k)
\label{Rough0}
\end{equation}%
and%
\begin{equation}
\mathbf{P}\left( 0\leq S_{n}<y,L_{n}\geq 0\right) \ \leq \
C\,b_{n}\,\int_{0}^{y}V(-z)dz\ .  \label{Rough22}
\end{equation}%
By these inequalities we deduce the estimates
\begin{eqnarray}
&&\mathbf{E}\left[ e^{S_{j}};S_{n}\leq \varphi (n);\tau _{n}=j\right]  \notag
\\
&=&\mathbf{E}\left[ e^{S_{j}}\mathbf{P}\left( S_{n-j}^{\prime }\leq \varphi
(n)-S_{j};L_{n-j}^{\prime }\geq 0|S_{j}\right) ;\tau _{j}=j\right]  \notag \\
&\leq &\sum_{k=1}^{\infty }e^{-k+1}\mathbf{P}\left( S_{j}\in \lbrack
-k,-k+1),\tau _{j}=j\right) \mathbf{P}\left( S_{n-j}\leq \varphi
(n)+k;L_{n-j}\geq 0\right)  \notag \\
&\leq &Cb_{j}\sum_{k=1}^{\infty }e^{-k+1}U(k)\mathbf{P}\left( S_{n-j}\leq
\varphi (n)+k;L_{n-j}\geq 0\right)  \label{RightTail0} \\
&\leq &Cb_{j}b_{n-j}\sum_{k=1}^{\infty }e^{-k}U(k)\int_{0}^{\varphi
(n)+k}V(-z)dz.\   \label{RightTail}
\end{eqnarray}%
According to Lemma \ref{L_estimV} there exists a constant $C(2)>0$ such that
\begin{equation*}
\int_{0}^{2\varphi (n)}V(-z)dz\leq C(2)\int_{0}^{\varphi (n)}V(-z)dz.
\end{equation*}%
Besides, $b_{j}\leq Cb_{n}$ for all $n/2\leq j\leq n$ in view of (\ref{Def_b}%
). Using the inequality $V(x+y)\leq V(x)+V(y)$ valid for $x\leq 0,y\leq 0$
(see, for instance, the proof of Corollary 2.4 in \cite{ABGV2011}) and (\ref%
{RightTail0}), we conclude that
\begin{eqnarray*}
&&\sum_{j\geq n/2}\mathbf{P}\left( Z_{n}>0;S_{n}\leq \varphi (n),\tau
_{n}=j\right) \\
&\leq &\sum_{j\geq n/2}\mathbf{E}\left[ e^{S_{j}},S_{n}\leq \varphi (n);\tau
_{n}=j\right] \\
&\leq &C\sum_{k=1}^{\infty }e^{-k}U(k)\sum_{j\geq n/2}b_{j}\mathbf{P}\left(
S_{n-j}\leq \varphi (n)+k;L_{n-j}\geq 0\right) \\
&\leq &Cb_{n}\sum_{k=1}^{\infty }e^{-k}U(k)V(-\varphi (n)-k) \\
&\leq &C_{1}b_{n}\left( \sum_{k=1}^{\infty }e^{-k}U(k)(V(-\varphi
(n))+V(-k))\right) \\
&\leq &C_{1}b_{n}\left( C_{2}V(-\varphi (n))+C_{3}\right) .
\end{eqnarray*}%
Since $V(-x),x>0,$ is a regularly varying function as $x\rightarrow \infty $
with a positive index, it follows by Corollary \ref{C_IntegVW} that
\begin{eqnarray}
&&\sum_{j\geq n/2}\mathbf{E}\left[ e^{S_{\tau _{n}}},S_{n}\leq \varphi
(n);\tau _{n}=j\right] \leq C_{3}b_{n}V(-\varphi (n))=o\left(
b_{n}\int_{0}^{\varphi (n)}V(-z)dz\right)  \notag \\
&=&o(\mathbf{P}\left( S_{n}\leq \varphi (n),L_{n}\geq 0\right) )
\label{Big_j}
\end{eqnarray}%
as $n\rightarrow \infty $. Further, in view of\ (\ref{RightTail})
\begin{eqnarray}
&&\sum_{J+1\leq j<n/2}\mathbf{E}\left[ e^{S_{j}},S_{n}\leq \varphi (n);\tau
_{n}=j\right]  \notag \\
&\leq &C\sum_{J\leq j<n/2}b_{j}b_{n-j}\sum_{k=1}^{\infty
}e^{-k}U(k)\int_{0}^{\varphi (n)+k}V(-z)dz  \notag \\
&\leq &Cb_{n}\sum_{J\leq j<n/2}b_{j}\left( \sum_{k=1}^{-\infty
}e^{-k}U(k)\left( \int_{0}^{2\varphi (n)}V(-z)dz+\int_{0}^{2k}V(-z)dz\right)
\right)  \notag \\
&\leq &Cb_{n}\sum_{J\leq j<n/2}b_{j}\left( C_{1}\int_{0}^{2\varphi
(n)}V(-z)dz+C_{2}\right) \leq C_{3}\sum_{J\leq j\leq \infty }b_{j}\times
b_{n}\int_{0}^{\varphi (n)}V(-z)dz  \notag \\
&=&\varepsilon _{J}\mathbf{P}\left( S_{n}\leq \varphi (n),L_{n}\geq 0\right)
,  \label{Small_j}
\end{eqnarray}%
where $\varepsilon _{J}\rightarrow 0$ as $J\rightarrow \infty $. Combining (%
\ref{Big_j}) and (\ref{Small_j}) and letting to infinity first $n$ and than $%
J$ we see that
\begin{equation*}
\lim_{J\rightarrow \infty }\lim_{n\rightarrow \infty }\frac{R(J,n)}{\mathbf{P%
}\left( S_{n}\leq \varphi (n),L_{n}\geq 0\right) }\leq C\lim_{J\rightarrow
\infty }\varepsilon _{J}=0.
\end{equation*}

Lemma \ref{L_remainder} is proved.

\begin{lemma}
\label{L_fixed_j}Under the conditions of Theorem \ref{T_smallDevi} for each
fixed $j$%
\begin{equation*}
\mathbf{P}\left( Z_{n}>0,S_{n}\leq \varphi (n),\tau _{n}=j\right) \sim
\Theta (j)g_{\alpha ,\beta }(0)b_{n}\int_{0}^{\varphi (n)}V(-z)dz
\end{equation*}%
as $n\rightarrow \infty ,$ where%
\begin{equation*}
\Theta (j)=\sum_{k=1}^{\infty }\mathbf{P}(Z_{j}=k,\tau _{j}=j)\mathbf{E}^{+}%
\left[ 1-F_{0,\infty }^{k}(0)\right] \leq \mathbf{P}(\tau _{j}=j)\leq 1.
\end{equation*}
\end{lemma}

\textbf{Proof}. First observe that for $0\leq j<n$ and any $\varepsilon \in
(0,1)$%
\begin{eqnarray*}
&&\mathbf{P}\left( Z_{n}>0,S_{n}\leq \varphi (n);\tau _{n}=j;S_{j}\leq
-\varepsilon \varphi (n)\right) \\
&=&\mathbf{E}\left[ 1-F_{0,n}(0),S_{n}\leq \varphi (n);\tau _{n}=j;S_{j}\leq
-\varepsilon \varphi (n)\right] \\
&\leq &\mathbf{E}\left[ e^{S_{\tau _{n}}},S_{n}\leq \varphi (n);\tau
_{n}=j;S_{j}\leq -\varepsilon \varphi (n)\right] \\
&=&\int_{-\infty }^{-\varepsilon \varphi (n)}e^{x}\mathbf{P}\left( S_{j}\in
dx,\tau _{j}=j\right) \mathbf{P}\left( S_{n-j}\leq \varphi (n)-x;L_{n-j}\geq
0\right) .
\end{eqnarray*}%
Using the decomposition $(-\infty ,-\varepsilon \varphi (n)]=(-\infty
,-\varphi (n)]\cup (-\varphi (n),-\varepsilon \varphi (n)]$ we have%
\begin{eqnarray*}
&&\int_{-\infty }^{-\varepsilon \varphi (n)}e^{x}\mathbf{P}\left( S_{j}\in
dx,\tau _{j}=j\right) \mathbf{P}\left( S_{n-j}\leq \varphi (n)-x;L_{n-j}\geq
0\right) \\
&\leq &\mathbf{P}\left( S_{n-j}\leq 2\varphi (n);L_{n-j}\geq 0\right)
\int_{-\varphi (n)}^{-\varepsilon \varphi (n)}e^{x}\mathbf{P}\left( S_{j}\in
dx,\tau _{j}=j\right) \\
&&+\int_{-\infty }^{-\varphi (n)}e^{x}\mathbf{P}\left( S_{j}\in dx,\tau
_{j}=j\right) \mathbf{P}\left( S_{n-j}\leq -2x;L_{n-j}\geq 0\right) .
\end{eqnarray*}%
Since $V(-x),x\geq 0,$ is a renewal function, $V(-x)\leq C(\left\vert
x\right\vert +1)$. Using this estimate and Lemma \ref{L_double} we conclude
that%
\begin{eqnarray*}
&&\int_{-\infty }^{-\varphi (n)}e^{x}\mathbf{P}\left( S_{j}\in dx,\tau
_{j}=j\right) \mathbf{P}\left( S_{n-j}\leq -2x;L_{n-j}\geq 0\right) \\
&\leq &Cb_{n-j}\int_{-\infty }^{-\varphi (n)}e^{x}\mathbf{P}\left( S_{j}\in
dx,\tau _{j}=j\right) \int_{0}^{-2x}V(-z)dz \\
&\leq &C_{1}b_{n-j}\int_{-\infty }^{-\varphi (n)}e^{x}\left\vert
x\right\vert ^{2}\mathbf{P}\left( S_{j}\in dx,\tau _{j}=j\right) \\
&\leq &C_{1}b_{n-j}\mathbf{E}\left[ e^{S_{j}}\left\vert S_{j}\right\vert
^{2}I\left\{ S_{j}\leq -\varphi (n)\right\} \right] =o(b_{n-j}).
\end{eqnarray*}%
This and the estimate for $V(-x)$ given above imply
\begin{eqnarray*}
&&\mathbf{P}\left( S_{n-j}\leq 2\varphi (n);L_{n-j}\geq 0\right)
\int_{-\varphi (n)}^{-\varepsilon \varphi (n)}e^{x}\mathbf{P}\left( S_{j}\in
dx,\tau _{j}=j\right) \\
&\leq &Cb_{n-j}\int_{0}^{2\varphi (n)}V(-z)dze^{-\varepsilon \varphi (n)} \\
&\leq &Cb_{n-j}\int_{0}^{2\varphi (n)}\left( z+1\right) dze^{-\varepsilon
\varphi (n)}=o(b_{n-j})
\end{eqnarray*}%
as $n\rightarrow \infty $. Therefore, for each fixed $j$ and $\varepsilon
\in (0,1)$%
\begin{equation*}
\mathbf{P}\left( Z_{n}>0,S_{n}\leq \varphi (n);\tau _{n}=j;S_{j}\leq
-\varepsilon \varphi (n)\right) =o(b_{n-j}).
\end{equation*}%
Further, we have
\begin{eqnarray*}
&&\mathbf{P}\left( Z_{n}>0,S_{n}\leq \varphi (n);\tau
_{n}=j,S_{j}>-\varepsilon \varphi (n)\right) \\
&=&\mathbf{E}\left[ \mathbf{P}\left( Z_{n}>0,S_{n}\leq \varphi (n);\tau
_{n}=j|S_{j},Z_{j}\right) ;S_{j}>-\varepsilon \varphi (n)\right] \\
&=&\sum_{k=1}^{\infty }\int_{-\varepsilon \varphi (n)}^{0}\mathbf{P}%
(Z_{j}=k,S_{j}\in dx,\tau _{j}=j)\mathbf{E}\left[ 1-F_{0,n-j}^{k}(0),S_{n-j}%
\leq \varphi (n)-x,L_{n-j}\geq 0\right] .
\end{eqnarray*}%
Given $\varphi (n)=o(a_{n})$ we see that for any $\varepsilon >0$ and all $%
x\in (-\varepsilon \varphi (n),0]$ and each $k\in \mathbf{N}$
\begin{eqnarray*}
&&\mathbf{E}\left[ 1-F_{0,n-j}^{k}(0),S_{n-j}\leq \varphi (n),L_{n-j}\geq 0%
\right] \\
&\leq &\mathbf{E}\left[ 1-F_{0,n-j}^{k}(0),S_{n-j}\leq \varphi
(n)-x,L_{n-j}\geq 0\right] \\
&\leq &\mathbf{E}\left[ 1-F_{0,n-j}^{k}(0),S_{n-j}\leq (1+\varepsilon
)\varphi (n),L_{n-j}\geq 0\right] .
\end{eqnarray*}%
Hence, using Lemma \ref{L_cond}, Corollary \ref{C_IntegVW} and the fact that
$\int_{0}^{x}V(-z)dz$ is a regularly varying function with index $\alpha
\rho +1$ we conclude that as $n\rightarrow \infty $%
\begin{equation*}
\lim_{n\rightarrow \infty }\inf_{x\in (-\varepsilon \varphi (n),0]}\frac{%
\mathbf{E}\left[ 1-F_{0,n-j}^{k}(0),S_{n-j}\leq \varphi (n)-x,L_{n-j}\geq 0%
\right] }{\mathbf{P}\left( S_{n-j}\leq \varphi (n),L_{n-j}\geq 0\right) }%
\geq \mathbf{E}^{+}\left[ 1-F_{0,\infty }^{k}(0)\right]
\end{equation*}%
and%
\begin{eqnarray*}
&&\lim_{n\rightarrow \infty }\sup_{x\in (-\varepsilon \varphi (n),0]}\frac{%
\mathbf{E}\left[ 1-F_{0,n-j}^{k}(0),S_{n-j}\leq \varphi (n)-x,L_{n-j}\geq 0%
\right] }{\mathbf{P}\left( S_{n-j}\leq \varphi (n),L_{n-j}\geq 0\right) } \\
&\leq &\left( 1+\varepsilon \right) ^{\alpha \rho +1}\mathbf{E}^{+}\left[
1-F_{0,\infty }^{k}(0)\right] .
\end{eqnarray*}

Therefore, for each fixed $k$ and $j$%
\begin{eqnarray*}
\lim_{\varepsilon \downarrow 0}\lim_{n\rightarrow \infty }\int_{-\varepsilon
\varphi (n)}^{0}\mathbf{P}(Z_{j} &=&k;S_{j}\in dx,\tau _{j}=j) \\
&&\times \frac{\mathbf{E}\left[ 1-F_{0,n-j}^{k}(0),S_{n-j}\leq \varphi
(n)-x,L_{n-j}\geq 0\right] }{\mathbf{P}\left( S_{n}\leq \varphi
(n),L_{n}\geq 0\right) } \\
&=&\int_{-\infty }^{0}\mathbf{P}(Z_{j}=k;S_{j}\in dx,\tau _{j}=j)\mathbf{E}%
^{+}\left[ 1-F_{0,\infty }^{k}(0)\right] \\
&=&\mathbf{P}(Z_{j}=k,\tau _{j}=j)\mathbf{E}^{+}\left[ 1-F_{0,\infty }^{k}(0)%
\right] .
\end{eqnarray*}%
Further,%
\begin{eqnarray*}
\sum_{k=K+1}^{\infty }\int_{-\varepsilon \varphi (n)}^{0}\mathbf{P}(Z_{j}
&=&k;S_{j}\in dx,\tau _{j}=j)\mathbf{E}\left[ 1-F_{0,n-j}^{k}(0),S_{n-j}\leq
\varphi (n)-x,L_{n-j}\geq 0\right] \\
&\leq &\sum_{k=K+1}^{\infty }\mathbf{P}(Z_{j}=k,\tau _{j}=j)\mathbf{P}\left(
S_{n-j}\leq (1+\varepsilon )\varphi (n),L_{n-j}\geq 0\right) \\
&\leq &\mathbf{P}(Z_{j}\geq K+1)\mathbf{P}\left( S_{n-j}\leq (1+\varepsilon
)\varphi (n),L_{n-j}\geq 0\right)
\end{eqnarray*}%
and, by Lemma \ref{L_estimV} the right-hand side of this relation is
\begin{equation*}
o\left( \mathbf{P}\left( S_{n}\leq \varphi (n),L_{n}\geq 0\right) \right)
\end{equation*}%
as $K\rightarrow \infty $. Hence, as $n\rightarrow \infty $
\begin{eqnarray*}
\mathbf{P}\left( Z_{n}>0,S_{n}\leq \varphi (n),\tau _{n}=j\right) &\sim
&\Theta (j)\mathbf{P}\left( S_{n}\leq \varphi (n),L_{n}\geq 0\right) \\
&\sim &\Theta (j)g_{\alpha ,\beta }(0)b_{n}\int_{0}^{\varphi (n)}V(-z)dz.
\end{eqnarray*}

for each fixed $j$.

Lemma \ref{L_fixed_j} is proved.

\textbf{Proof of Theorem \ref{T_smallDevi}}. Combining Lemmas \ref%
{L_remainder} and \ \ref{L_fixed_j} we see that
\begin{equation*}
\mathbf{P}\left( Z_{n}>0,S_{n}\leq \varphi (n)\right) \sim \Theta g_{\alpha
,\beta }(0)b_{n}\int_{0}^{\varphi (n)}V(-z)dz,
\end{equation*}%
where%
\begin{equation}
\Theta =\sum_{j=0}^{\infty }\Theta (j)==\sum_{j=1}^{\infty
}\sum_{k=1}^{\infty }\mathbf{P}(Z_{j}=k,\tau _{j}=j)\mathbf{E}^{+}\left[
1-F_{0,\infty }^{k}(0)\right] .  \label{Def_Theta}
\end{equation}

Positivity and finiteness of $\Theta $ was proved in \cite{agkv}, Theorem
1.1 (see, in particular, formula (4.10)).

Theorem \ref{T_smallDevi} is proved.

\end{document}